\documentclass{ws-jaa}
\usepackage{amssymb}    % AMS fonts
\usepackage{enumitem}
\usepackage{rotating}
\usepackage{url}

\numberwithin{equation}{section}

\newcommand{\seclabel}[1]{\label{sec:#1}}   % section
\newcommand{\thmlabel}[1]{\label{thm:#1}}   % theorem
\newcommand{\lemlabel}[1]{\label{lem:#1}}   % lemma
\newcommand{\corlabel}[1]{\label{cor:#1}}   % corollary
\newcommand{\prplabel}[1]{\label{prp:#1}}   % proposition
\newcommand{\deflabel}[1]{\label{def:#1}}   % definition
\newcommand{\secref}[1]{\ref{sec:#1}}   % section
\newcommand{\thmref}[1]{\ref{thm:#1}}   % theorem
\newcommand{\lemref}[1]{\ref{lem:#1}}   % lemma
\newcommand{\corref}[1]{\ref{cor:#1}}   % corollary
\newcommand{\prpref}[1]{\ref{prp:#1}}   % proposition
\newcommand{\inv}{^{-1}}
\newcommand{\Bol}{\textsc{Bol}}
\newcommand{\Mlt}{\mathrm{Mlt}}
\newcommand{\Inn}{\mathrm{Inn}}
\newcommand{\LMlt}{\mathrm{LMlt}}
\newcommand{\RMlt}{\mathrm{RMlt}}
\newcommand{\LInn}{\mathrm{LInn}}
\newcommand{\RInn}{\mathrm{RInn}}

\newcommand{\ZZ}{\mathbb{Z}}
\newcommand{\KK}{\mathcal{K}}
\newcommand{\HH}{\mathcal{H}}
\newcommand{\TT}{\mathcal{T}}
\newcommand{\BT}{\mathcal{BT}}
\newcommand{\PP}{\mathcal{P}}
\newcommand{\Aut}{\mathrm{Aut}}

\catchline{}{}{}{}{}
\begin{document}

\title{Uniquely $2$-divisible Bol loops}

\author{Tuval Foguel}

\address{Department of Mathematics and Computer Science \\
Stillwell 426 \\
Western Carolina University \\
Cullowhee, NC 28723 USA\\
\email{tsfoguel@wcu.edu}}

\author{Michael Kinyon}

\address{Department of Mathematics\\
2360 S Gaylord St \\
University of Denver\\
Denver CO 80208 USA \\
\email{mkinyon@math.du.edu}}

\maketitle

\begin{history}
\received{(Day Month Year)}
\revised{(Day Month Year)}
\accepted{(Day Month Year)}
\comby{(xxxxxxxxx)}
\end{history}

\begin{abstract}
Although any finite Bol loop of odd prime exponent is
solvable, we show there exist such Bol loops with trivial center.
We also construct finitely generated, infinite, simple Bruck loops of odd
prime exponent for sufficiently large primes. This shows that the
Burnside problem for Bruck loops has a negative answer.
\end{abstract}


\ccode{2000 Mathematics Subject Classification: 20N05}

\section{Introduction}
\seclabel{intro}

A \emph{loop} $(Q,\cdot)$ consists of a set $Q$ with a binary operation $\cdot : Q\times Q\to Q$
such that (i) for all $a,b\in Q$, the equations $ax = b$ and $ya = b$ have unique solutions
$x,y\in Q$, and (ii) there exists $1\in Q$ such that $1x = x1 = x$ for all $x\in Q$.
A loop $Q$ is said to be a (right) \emph{Bol loop} if it satisfies the identity
\[
[(xy)z]y = x[(yz)y]  \tag{\Bol}
\]
for all $x,y,z\in Q$. For $x\in Q$, define the \emph{right} and \emph{left translations} by $x$
by, respectively, $y R_x = yx$ and $yL_x = xy$ for all $y\in Q$. Then the Bol identity is
equivalent to
\[
R_x R_y R_x = R_{(xy)x}
\]
for all $x,y\in Q$. Much of the literature on Bol loops (\emph{e.g.}, \cite{Kiechle}) works with the dual notion of
\emph{left} Bol loop defined by the identity $x[y(xz)] = [x(yx)]z$. Results for left Bol
loops dualize trivially to results about right Bol loops.
For background in loop theory, we refer the reader to \cite{Br, Pf}.
Basic information about Bol loops, especially facts not explicitly cited below, can
be found in \cite{Rob66}.

Let $Q$ be a loop. For $x\in Q$ and $n\in \ZZ$, set $x^n = 1 R_x^n$.
If $Q$ is a Bol loop, then $Q$ is \emph{power-associative}, that is,
$x^m x^n = x^{m+n}$ for all $x\in Q$, $m,n\in \ZZ$, and $Q$ is also
\emph{right power-alternative}, that is, $R_x^n = R_{x^n}$ for all $x\in Q$,
$n\in \ZZ$. Finite Bol loops satisfy the \emph{elementwise Lagrange property},
that is, the order of any element divides the order of the loop. (It is still
an open whether or
not Bol loops satisfy the full Lagrange property: does the order of a subloop
of a finite Bol loop divide the order of the loop?)
A Bol loop is said to be \emph{uniquely} $2$-\emph{divisible} if the
squaring map $x\mapsto x^2$ is a bijection. If $Q$ is a finite Bol
loop, unique $2$-divisibility is equivalent to every element of $Q$ having
finite odd order, and also to $Q$ itself having odd order \cite{FKP}.
Bol loops of odd order satisfy the \emph{Cauchy property}, that is, for
each prime $p$, there exists an element of order $p$ \cite{FKP}.

For $p$ a prime,
a power-associative loop is called a $p$-\emph{loop} if every element
has finite order which is some power of $p$, that is, for each $x$, there
exists $n > 0$ such that $x^{p^n} = 1$. If $Q$ is a finite Bol loop, then
by the elementwise Lagrange and Cauchy properties, being a $p$-loop is equivalent
to $Q$ itself having order some power of $p$.

Two important subclasses of Bol loops are \emph{Moufang loops} and
\emph{Bruck loops}. The former can be defined in various ways, for instance,
as those Bol loops satisfying the antiautomorphic inverse property
$(xy)\inv = y\inv x\inv$. Bruck loops are Bol loops satisfying the
automorphic inverse property $(xy)\inv = x\inv y\inv$.
The intersection of these two classes is the class of \emph{commutative
Moufang loops}.

The \emph{multiplication group} of a loop $Q$ is the group
$\Mlt(Q) = \langle R_x, L_x | x \in Q\rangle$. The \emph{inner mapping group}
$\Inn(Q)$ is the stabilizer of the identity element $1\in Q$.
A subloop which is invariant under the action of $\Inn(Q)$
is said to be \emph{normal}. An important fact about normal subloops
we will use later is that they are \emph{permutable}, that is,
if $H, K$ are subloops of a loop $Q$ and if $K$ is normal, then
$HK$ and $KH$ are subloops of $Q$. A loop with no nontrivial
normal subloops is \emph{simple}.

The \emph{right nucleus} of a loop $Q$ is the set
$N_r(Q) = \{ a\in Q\ |\ (xy)a = x(ya)\ \forall x,y\in Q\}$,
the \emph{left nucleus} is
$N_l(Q) = \{ a\in Q\ |\ (ax)y = a(xy)\ \forall x,y\in Q\}$,
the \emph{middle nucleus} is
$N_m(Q) = \{ a\in Q\ |\ (xa)y = x(ay)\ \forall x,y\in Q\}$,
and the \emph{nucleus} is $N(Q) = N_r(Q)\cap N_m(Q)\cap N_l(Q)$.
These are all subloops. In a Bol loop $Q$, $N_r(Q) = N_m(Q)$
is a normal subloop \cite{FKP}, but $N_l(Q)$ is not necessarily normal.

For a subset $S$ of a loop $Q$, we borrow terminology from semigroup
theory and define the \emph{commutant of} $S$ \emph{in} $Q$ to be the set
$C_Q(S) = \{ a\in Q \ |\ ax = xa\ \forall x\in S\}$. As usual, if $S =\{x\}$
is a singleton, we write $C_Q(x)$ instead of $C_Q(\{x\})$.
The \emph{commutant} of $Q$ itself is $C(Q) = C_Q(Q)$.
(This latter set is also
known as the centrum or semicenter or commutative center.)
Note that $C(Q) = \bigcap_{x\in Q} C_Q(x)$. The commutant
of a Moufang loop is a subloop \cite{Pf}, but this is
not generally true in Bol loops or even Bruck loops \cite{KPV}. It is, however,
true that in a uniquely $2$-divisible Bol loop, the commutant is
a subloop \cite{KP}. Even in Moufang loops, the normality of
the commutant remains an open problem. In \S\secref{comm}, we
show that \emph{in a uniquely} $2$-\emph{divisible Bruck loop, the commutant is
a normal subloop.}

The \emph{center} of a loop $Q$ is $Z(Q) = N(Q)\cap C(Q)$. This is
exactly the fixed point subset of $\Inn(Q)$.
The notion of \emph{upper central series} is defined in the
same way as in group theory, and a loop is \emph{centrally nilpotent}
if it has a finitely terminating upper central series.
As in group theory, a loop is \emph{solvable} if it has a normal series
with factors which are abelian groups. Centrally nilpotent loops are
necessarily solvable.

Much has been learned about Bol loops in the last few years,
especially thanks to work of Nagy. For about three decades,
the existence of finite, simple, non-Moufang Bol loops was the
main open problem in loop theory. Firstly, Nagy gave a very
general construction of a large class of such loops \cite{Nagy0}. The
construction showed that not only are there many such loops
(more than will likely allow for a classification any time soon),
but there are also surprising examples. For instance, Nagy showed
that there exist finite, simple, non-Moufang Bol loops of odd
order \cite{Nagy0}. This answered negatively a question of \cite{FKP},
and is in sharp contrast to the more specialized cases: every finite
Bruck loop of odd order is solvable \cite{Gl1} and every finite
Moufang loop of odd order is solvable \cite{Gl2}.

Turning now to $p$-loops, we consider first $p = 2$. Finite
Moufang $2$-loops are centrally nilpotent \cite{GW}.
There exist finite Bruck loops of exponent $2$ with trivial
center \cite{KieN}; the smallest order at which this occurs is $16$.
More recently, building on earlier work of
Aschbacher \cite{Asch2}, Nagy constructed a finite simple Bol loop of
exponent $2$ and order $96$ \cite{Nagy2}.

What of Bol $p$-loops for $p$ an odd prime?
Thanks to seminal work of Glauberman, it is known
that every finite Bruck $p$-loop has nontrivial center \cite{Gl1},
and every finite Moufang $p$-loop has nontrivial center \cite{Gl2}.
It follows immediately by induction that such loops in either class
are centrally nilpotent.

In the general case, the best positive result for a finite Bol $p$-loop is
again due to Nagy: such a loop is solvable \cite{Nagy1}. However, until recently,
the following problem was considered to be open:
\begin{center}
\emph{Does every finite Bol} $p$-\emph{loop have nontrivial center?}
\end{center}
\noindent In $\S\secref{2}$ we present examples for $p = 3$ showing
that the answer to this question is no. As it turns out, our examples
are related to a loop which has been in the literature since 1963, and
which, suitably interpreted, would have already answered the question.

Finitely generated commutative Moufang loops are finite \cite{Br},
and so it follows from Glauberman's results that for each odd prime $p$,
a finitely generated
commutative Moufang $p$-loop of finite exponent is centrally nilpotent.
(This can, in fact, be shown directly for commutative Moufang loops \cite{Br}.)
One might then ask  the following:
\begin{center}
\emph{Is every finitely generated Bruck} $p$-\emph{loop of finite exponent
centrally nilpotent?}
\end{center}
\noindent In \S\secref{3}, we answer this question negatively by constructing
a family of Bruck loops we call \emph{Bruck-Tarski monsters}, which are
finitely generated, infinite and simple.

\section{Finite Bol $3$-loops with trivial center}
\seclabel{2}

While searching for particular types of projective planes,
Keedwell gave an example of a normalized latin square of order $27$ with
certain properties. He published the example in two papers
\cite{Keed0, Keed1} in 1963 and 1965, respectively.
D\'{e}nes and Keedwell later included the example in their
well-known book on latin squares (\cite{DK}, p. 50).

A normalized latin square is the Cayley table for a loop.
It is obvious by inspection that the loop from Keedwell's
latin square has exponent $3$, but no other properties are immediately
apparent. Perhaps because the example appeared in the latin square
literature, the loop theory community seems to have been unaware of it.

Using the LOOPS package \cite{NV} for GAP \cite{GAP}, we analyzed
Keedwell's loop and found that it was, in fact, a \emph{left} Bol
loop with trivial center. The loop has a normal associative
subloop of order $9$. Using tools from the LOOPS package,
we took the transpose of the Cayley table to get a \emph{right}
Bol loop, and then constructed an isomorphic copy in which the
normal subloop of order $9$ is in the upper left corner.

Two loop structures $(Q,\cdot)$ and $(Q,\circ)$ on the same
underlying set $Q$ are said to be (principally) \emph{isotopic} if
there exist $a,b\in Q$ such that $x\circ y =
xR_a\inv\cdot yL_b\inv$ where the translations are taken in
$(Q,\cdot)$. Every isotope of a Bol loop is itself a Bol loop.
Further, every loop isotopic to a Bol loop $(Q,\cdot)$ is
isomorphic to an isotope of a particular form, namely one given by
$x\circ y = xR_a \cdot yL_a\inv$ for $a\in Q$ \cite{Rob66}.

We computed all principal isotopes of our modified version of
Keedwell's loop, a task made
easier since it is enough to consider isotopes of the form
$x\circ y = xR_a \cdot yL_a\inv$. It turns out that
there are exactly two isotopy classes. We present here a
representative from each class in Tables 1 and 2.

\begin{sidewaystable}
\begin{center}{\small
\[
\begin{tabular}{ccccccccccccccccccccccccccc}
\\ \\ \\ \\ \\ \\ \\ \\ \\ \\ \\ \\ \\ \\ \\ \\ \\ \\ \\ \\ \\ \\ \\ \\ \\ \\ \\ \\ \\ \\
1& 2&  3&  4&  5&  6&  7&  8&  9& 10& 11& 12& 13& 14& 15& 16& 17& 18& 19& 20& 21& 22& 23& 24& 25& 26& 27
 \\  2&  3&  1&  5&  6&  4&  8&  9&  7& 15& 13& 14& 18& 16& 17& 12& 10& 11& 22& 23& 24& 25& 26& 27& 19& 20& 21
 \\  3&  1&  2&  6&  4&  5&  9&  7&  8& 17& 18& 16& 11& 12& 10& 14& 15& 13& 25& 26& 27& 19& 20& 21& 22& 23& 24
 \\  4&  5&  6&  7&  8&  9&  1&  2&  3& 14& 15& 13& 17& 18& 16& 11& 12& 10& 26& 27& 25& 20& 21& 19& 23& 24& 22
 \\  5&  6&  4&  8&  9&  7&  2&  3&  1& 16& 17& 18& 10& 11& 12& 13& 14& 15& 20& 21& 19& 23& 24& 22& 26& 27& 25
 \\  6&  4&  5&  9&  7&  8&  3&  1&  2& 12& 10& 11& 15& 13& 14& 18& 16& 17& 23& 24& 22& 26& 27& 25& 20& 21& 19
 \\  7&  8&  9&  1&  2&  3&  4&  5&  6& 18& 16& 17& 12& 10& 11& 15& 13& 14& 24& 22& 23& 27& 25& 26& 21& 19& 20
 \\  8&  9&  7&  2&  3&  1&  5&  6&  4& 11& 12& 10& 14& 15& 13& 17& 18& 16& 27& 25& 26& 21& 19& 20& 24& 22& 23
 \\  9&  7&  8&  3&  1&  2&  6&  4&  5& 13& 14& 15& 16& 17& 18& 10& 11& 12& 21& 19& 20& 24& 22& 23& 27& 25& 26
 \\ 10& 11& 12& 13& 14& 15& 16& 17& 18& 19& 21& 20& 25& 27& 26& 22& 24& 23&  1&  3&  2&  7&  9&  8&  4&  6&  5
 \\ 11& 12& 10& 14& 15& 13& 17& 18& 16& 27& 26& 25& 24& 23& 22& 21& 20& 19&  8&  7&  9&  5&  4&  6&  2&  1&  3
 \\ 12& 10& 11& 15& 13& 14& 18& 16& 17& 23& 22& 24& 20& 19& 21& 26& 25& 27&  6&  5&  4&  3&  2&  1&  9&  8&  7
 \\ 13& 14& 15& 16& 17& 18& 10& 11& 12& 21& 20& 19& 27& 26& 25& 24& 23& 22&  9&  8&  7&  6&  5&  4&  3&  2&  1
 \\ 14& 15& 13& 17& 18& 16& 11& 12& 10& 26& 25& 27& 23& 22& 24& 20& 19& 21&  4&  6&  5&  1&  3&  2&  7&  9&  8
 \\ 15& 13& 14& 18& 16& 17& 12& 10& 11& 22& 24& 23& 19& 21& 20& 25& 27& 26&  2&  1&  3&  8&  7&  9&  5&  4&  6
 \\ 16& 17& 18& 10& 11& 12& 13& 14& 15& 20& 19& 21& 26& 25& 27& 23& 22& 24&  5&  4&  6&  2&  1&  3&  8&  7&  9
 \\ 17& 18& 16& 11& 12& 10& 14& 15& 13& 25& 27& 26& 22& 24& 23& 19& 21& 20&  3&  2&  1&  9&  8&  7&  6&  5&  4
 \\ 18& 16& 17& 12& 10& 11& 15& 13& 14& 24& 23& 22& 21& 20& 19& 27& 26& 25&  7&  9&  8&  4&  6&  5&  1&  3&  2
 \\ 19& 20& 21& 22& 23& 24& 25& 26& 27&  1&  7&  4&  6&  3&  9&  8&  5&  2& 10& 18& 14& 12& 17& 13& 11& 16& 15
 \\ 20& 21& 19& 23& 24& 22& 26& 27& 25&  5&  2&  8&  7&  4&  1&  3&  9&  6& 16& 15& 11& 18& 14& 10& 17& 13& 12
 \\ 21& 19& 20& 24& 22& 23& 27& 25& 26&  9&  6&  3&  2&  8&  5&  4&  1&  7& 13& 12& 17& 15& 11& 16& 14& 10& 18
 \\ 22& 23& 24& 25& 26& 27& 19& 20& 21&  2&  8&  5&  4&  1&  7&  9&  6&  3& 15& 11& 16& 14& 10& 18& 13& 12& 17
 \\ 23& 24& 22& 26& 27& 25& 20& 21& 19&  6&  3&  9&  8&  5&  2&  1&  7&  4& 12& 17& 13& 11& 16& 15& 10& 18& 14
 \\ 24& 22& 23& 27& 25& 26& 21& 19& 20&  7&  4&  1&  3&  9&  6&  5&  2&  8& 18& 14& 10& 17& 13& 12& 16& 15& 11
 \\ 25& 26& 27& 19& 20& 21& 22& 23& 24&  3&  9&  6&  5&  2&  8&  7&  4&  1& 17& 13& 12& 16& 15& 11& 18& 14& 10
 \\ 26& 27& 25& 20& 21& 19& 23& 24& 22&  4&  1&  7&  9&  6&  3&  2&  8&  5& 14& 10& 18& 13& 12& 17& 15& 11& 16
 \\ 27& 25& 26& 21& 19& 20& 24& 22& 23&  8&  5&  2&  1&  7&  4&  6&  3&  9& 11& 16& 15& 10& 18& 14& 12& 17& 13
\end{tabular}
\]}
\caption{Right Bol loop $1$ of order $27$ and exponent $3$ with trivial center.}
\end{center}
\end{sidewaystable}

\begin{sidewaystable}
\begin{center}{\small
\[
\begin{tabular}{ccccccccccccccccccccccccccc}
\\ \\ \\ \\ \\ \\ \\ \\ \\ \\ \\ \\ \\ \\ \\ \\ \\ \\ \\ \\ \\ \\ \\ \\ \\ \\ \\ \\ \\ \\
\\ 1&  2&  3&  4&  5&  6&  7&  8&  9& 10& 11& 12& 13& 14& 15& 16& 17& 18& 19& 20& 21& 22& 23& 24& 25& 26& 27
 \\  2&  3&  1&  5&  6&  4&  8&  9&  7& 15& 13& 14& 18& 16& 17& 12& 10& 11& 22& 23& 24& 25& 26& 27& 19& 20& 21
 \\  3&  1&  2&  6&  4&  5&  9&  7&  8& 17& 18& 16& 11& 12& 10& 14& 15& 13& 25& 26& 27& 19& 20& 21& 22& 23& 24
 \\  4&  5&  6&  7&  8&  9&  1&  2&  3& 14& 15& 13& 17& 18& 16& 11& 12& 10& 26& 27& 25& 20& 21& 19& 23& 24& 22
 \\  5&  6&  4&  8&  9&  7&  2&  3&  1& 16& 17& 18& 10& 11& 12& 13& 14& 15& 20& 21& 19& 23& 24& 22& 26& 27& 25
 \\  6&  4&  5&  9&  7&  8&  3&  1&  2& 12& 10& 11& 15& 13& 14& 18& 16& 17& 23& 24& 22& 26& 27& 25& 20& 21& 19
 \\  7&  8&  9&  1&  2&  3&  4&  5&  6& 18& 16& 17& 12& 10& 11& 15& 13& 14& 24& 22& 23& 27& 25& 26& 21& 19& 20
 \\  8&  9&  7&  2&  3&  1&  5&  6&  4& 11& 12& 10& 14& 15& 13& 17& 18& 16& 27& 25& 26& 21& 19& 20& 24& 22& 23
 \\  9&  7&  8&  3&  1&  2&  6&  4&  5& 13& 14& 15& 16& 17& 18& 10& 11& 12& 21& 19& 20& 24& 22& 23& 27& 25& 26
 \\ 10& 14& 18& 11& 15& 16& 12& 13& 17& 19& 20& 21& 22& 23& 24& 25& 26& 27&  1&  4&  7&  8&  2&  5&  6&  9&  3
 \\ 11& 15& 16& 12& 13& 17& 10& 14& 18& 27& 25& 26& 21& 19& 20& 24& 22& 23&  8&  2&  5&  6&  9&  3&  1&  4&  7
 \\ 12& 13& 17& 10& 14& 18& 11& 15& 16& 23& 24& 22& 26& 27& 25& 20& 21& 19&  6&  9&  3&  1&  4&  7&  8&  2&  5
 \\ 13& 17& 12& 14& 18& 10& 15& 16& 11& 21& 19& 20& 24& 22& 23& 27& 25& 26&  9&  3&  6&  4&  7&  1&  2&  5&  8
 \\ 14& 18& 10& 15& 16& 11& 13& 17& 12& 26& 27& 25& 20& 21& 19& 23& 24& 22&  4&  7&  1&  2&  5&  8&  9&  3&  6
 \\ 15& 16& 11& 13& 17& 12& 14& 18& 10& 22& 23& 24& 25& 26& 27& 19& 20& 21&  2&  5&  8&  9&  3&  6&  4&  7&  1
 \\ 16& 11& 15& 17& 12& 13& 18& 10& 14& 20& 21& 19& 23& 24& 22& 26& 27& 25&  5&  8&  2&  3&  6&  9&  7&  1&  4
 \\ 17& 12& 13& 18& 10& 14& 16& 11& 15& 25& 26& 27& 19& 20& 21& 22& 23& 24&  3&  6&  9&  7&  1&  4&  5&  8&  2
 \\ 18& 10& 14& 16& 11& 15& 17& 12& 13& 24& 22& 23& 27& 25& 26& 21& 19& 20&  7&  1&  4&  5&  8&  2&  3&  6&  9
 \\ 19& 21& 20& 25& 27& 26& 22& 24& 23&  1&  2&  3&  4&  5&  6&  7&  8&  9& 10& 12& 11& 16& 18& 17& 13& 15& 14
 \\ 20& 19& 21& 26& 25& 27& 23& 22& 24&  5&  6&  4&  8&  9&  7&  2&  3&  1& 16& 18& 17& 13& 15& 14& 10& 12& 11
 \\ 21& 20& 19& 27& 26& 25& 24& 23& 22&  9&  7&  8&  3&  1&  2&  6&  4&  5& 13& 15& 14& 10& 12& 11& 16& 18& 17
 \\ 22& 24& 23& 19& 21& 20& 25& 27& 26&  2&  3&  1&  5&  6&  4&  8&  9&  7& 15& 14& 13& 12& 11& 10& 18& 17& 16
 \\ 23& 22& 24& 20& 19& 21& 26& 25& 27&  6&  4&  5&  9&  7&  8&  3&  1&  2& 12& 11& 10& 18& 17& 16& 15& 14& 13
 \\ 24& 23& 22& 21& 20& 19& 27& 26& 25&  7&  8&  9&  1&  2&  3&  4&  5&  6& 18& 17& 16& 15& 14& 13& 12& 11& 10
 \\ 25& 27& 26& 22& 24& 23& 19& 21& 20&  3&  1&  2&  6&  4&  5&  9&  7&  8& 17& 16& 18& 14& 13& 15& 11& 10& 12
 \\ 26& 25& 27& 23& 22& 24& 20& 19& 21&  4&  5&  6&  7&  8&  9&  1&  2&  3& 14& 13& 15& 11& 10& 12& 17& 16& 18
 \\ 27& 26& 25& 24& 23& 22& 21& 20& 19&  8&  9&  7&  2&  3&  1&  5&  6&  4& 11& 10& 12& 17& 16& 18& 14& 13& 15 \end{tabular}
\]}
 \caption{Right Bol loop $2$ of order $27$ and exponent $3$ with trivial center.}
 \end{center}
\end{sidewaystable}

We computed various structures attached to these loops. Besides the multiplication
groups and inner mappings groups, we also computed the following.
The \emph{right} and \emph{left multiplication
groups} of a loop $Q$ are $\RMlt(Q) = \langle R_x \ |\ x\in Q\rangle$
and $\LMlt(Q) = \langle L_x\ |\ x\in Q\rangle$, and their respective
stabilizers of the identity element are the \emph{right} and \emph{left
inner mapping groups} $\LInn(Q)$ and $\RInn(Q)$.
Isotopic loops have isomorphic multiplication groups (both full and one-sided).
The \emph{associator subloop} $A(Q)$ is the smallest normal subloop
of $Q$ such that $Q/A(Q)$ is a group. The \emph{derived subloop}
$Q'$ (sometimes called the commutator-associator subloop) is the
smallest normal subloop of $Q$ such that $Q/Q'$ is an abelian group.

In the discussion that follows, $Q$ refers to either of the two loops.

\begin{itemize}
\item $G := \Mlt(Q) = \LMlt(Q)$ is a solvable group of order $2^6 3^7$ and derived
length $3$. The derived subgroup $G'$ has order $2^6 3^6$, while $G''$ is
elementary abelian of order $3^5$.
$G$ itself is a semidirect product of $C_3$ with $G'$ and is also a semidirect product of
a subgroup of order $2^6 3$ (the group [192,3], according to GAP \cite{GAP}) with $G''$.

(Note that if the order of $G$ were a power of $3$, then by
a theorem of A. A. Albert, the centers of the two loops
would be nontrivial \cite{Alb}.)

\item $H:= \Inn(Q) = \LInn(Q)$ is a solvable group
of order $2^6 3^4$ and derived length $2$. $H'$
is elementary abelian of order $3^4$ and $H$ itself factors as a semidirect
product of an abelian subgroup of order $2^6$ with $H'$.

\item $K:= \RMlt(Q)$ is a nilpotent group of order $3^5$ and nilpotency class $2$.
GAP \cite{GAP} identifies it as group $[243,37]$, a group with rank $3$.
$K' = Z(K)$ is elementary abelian of order $9$ and $K/Z(K)$ is also elementary abelian.

\item $\RInn(Q)$ is an elementary abelian group of order $9$.

\item $N_r(Q) = N_m(Q) = 1$.

\item $C(Q) = 1$.

\item $N_l(Q) = A(Q) = Q'$ is an elementary abelian group of order $9$
(the upper left corner of both tables). This is the unique nontrivial normal subloop.
\end{itemize}

\section{Commutants of uniquely $2$-divisible Bruck loops}
\seclabel{comm}

As mentioned in the introduction, in a uniquely $2$-divisible Bol loop $Q$,
the commutant $C(Q)$ is a subloop. In this section we address the normality
of this subloop in uniquely $2$-divisible Bruck loops.

A subloop is normal if it is invariant under the action of $\Inn(Q)$.
The group $\Inn(Q)$ is generated by the following permutations:
\[
R_{x,y} = R_x R_y R_{xy}\inv
\qquad
T_x = R_x L_x\inv
\qquad
L_{x,y} = L_x L_y L_{yx}\inv
\]
for all $x,y\in Q$. In a Bruck loop, the right inner mappings
$R_{x,y}$ play a special role. The following result can be found
in various sources; see, \emph{e.g.}, \cite{Kiechle} (translating
between left and right Bruck loops).

\begin{proposition}
In a Bruck loop, each $R_{x,y}$ is an automorphism.
\end{proposition}

\begin{theorem}
\thmlabel{comm-norm}
Let $Q$ be a uniquely $2$-divisible Bol loop such that each
$R_{x,y}$ is an automorphism. Then $C(Q)$ is a normal subloop.
\end{theorem}

\begin{proof}
The commutant $C(Q)$ is fixed by each $T_x$, $x\in Q$.
It is clear from the definition of $C(Q)$ that it is
invariant under the action of automorphisms, and so it follows that
$c R_{x,y}\in C(Q)$ for all $x,y\in Q$, $c\in C(Q)$.
It is thus enough to check that for all $c\in C(Q)$, $x,y\in Q$,
we have $c L_{x,y} \in C(Q)$. Now for fixed $x,y\in Q$, set
$z = c L_{x,y}$, and let $u = y^{1/2}$. We compute
\[
(xy)z = x(yc) = x(cy) = x[(cu)u] = x[(uc)u]
[(xu)c]u = [c(xu)]u\,
\]
using the right alternative property in the third equality and
(\Bol) in the fifth equality. Thus
\[
[(xy)z](xy)\inv = \{[c(xu)]u\}(xy)\inv = \{[c(xu)]u\}(xu^2)\inv =
c R_{xu,u}\,.
\]
Therefore $[(xy)z](xy)\inv \in C(Q)$. But then
$(xy)z = \{[(xy)z](xy)\inv\} (xy) = (xy)\{[(xy)z](xy)\inv\}$.
Canceling, we have $z = [(xy)z](xy)\inv$, and so $z\in C(Q)$, as claimed.
\end{proof}

\begin{corollary}
\corlabel{bruck-comm}
If $Q$ is a uniquely $2$-divisible Bruck loop, then $C(Q)$ is a normal subloop.
\end{corollary}

\section{Twisted subgroups and uniquely $2$-divisible Bruck loops}

In this section, we briefly review some facts about twisted subgroups,
and a standard construction of Bruck loops.

Recall that a subset $T$ of a group $G$ is a \textit{twisted subgroup} of $G$ if
(i) $1 \in T$, (ii) $T\inv = T$, and (iii) $xTx\subseteq T$ for all $x \in T$
\cite{Asch,FKP}.

Let $G$ be a group and suppose $\tau\in \Aut(G)$ has order $2$.
A typical example of a twisted subgroup which is not a subgroup is the set
$K(\tau) = \{ g\in G\ |\ g^{\tau} = g\inv \}$.

If $T$ is a twisted subgroup, the
\emph{Aschbacher radical} of $T$ is the set
$T' = \{ x_1 \cdots x_n \ |\ x_1\inv \cdots x_n\inv = 1, x_i\in T \}$.
This is a normal subgroup of $G$ \cite{Asch}, and it also measures whether or not $T$
arises as a set of inverted elements of an involutory automorphism in the
case where $T$ generates $G$.

\begin{proposition}[\cite{Asch}, Thm. 2.2; \cite{FKP}, Prop. 3.9]
\prplabel{rad}
Let $G$ be a group, $T\subseteq G$ a twisted subgroup and assume $\langle T\rangle = G$.
There exists $\tau\in \Aut(G)$ with $\tau^2 = 1$ such that $T\subseteq K(\tau)$ if and only
if $T' = 1$.
\end{proposition}

A twisted subgroup $T$ of a group $G$ is said to be \emph{uniquely}
$2$-\emph{divisible} if each $x\in T$ has a unique square root in $T$, that is,
a unique element $x^{1/2}\in T$ such that $(x^{1/2})^2 = x$. On such a
twisted subgroup, we may
define a new binary operation $\odot : T\times T\to T$ by
\[
x\odot y := (y x^2 y)^{1/2}
\]
for $x,y\in T$. Following Glauberman \cite{Gl1}, we denote the
magma $(T,\odot)$ by $T(1/2)$. The following properties are straightforward
to verify \cite{Gl1,FKP}. (Note that in \cite{Gl1, FKP}, the operation
is $x\odot y = (xy^2x)^{1/2}$, which gives a \emph{left} Bruck loop.)

\begin{itemize}
\item $T(1/2)$ is a uniquely $2$-divisible Bruck loop.
\item Integer powers of elements in $T$ formed in $G$ agree with those in $T(1/2)$.
Thus an element has finite order in $T$ if and only if it has the same order in $T(1/2)$.
\item If $K\subseteq T$ is a subloop of $T(1/2)$, then $K$ is a twisted subgroup of $G$.
\end{itemize}

\section{Bruck-Tarski Monsters}
\seclabel{3}

A \emph{Tarski monster} $\TT$ is an infinite group in which every nontrivial proper subgroup
is cyclic of order a fixed prime $p$. A Tarski monster is a simple group. Such groups were
constructed by Ol'shanskii \cite{Ol1} for every $p > 10^{75}$.
For more detail about these groups, see \cite{Ol2}.

One can clearly make sense of the notion of a Tarski monster for any class of
power-associative loops. Here we consider the following.

\begin{definition}
\deflabel{bruck-tarski}
A \emph{Bruck-Tarski monster} is an infinite Bruck loop in which every nontrivial
proper subloop is cyclic of prime order $p$.
\end{definition}

An immediate consequence of this definition is
that \emph{a Bruck-Tarski monster is a two-generated loop of exponent $p$}.

\begin{theorem}
\thmlabel{bt-simple}
A Bruck-Tarski monster $\BT$ is a simple loop.
\end{theorem}

\begin{proof}
Let $\KK, \HH$ be distinct proper subloops of $\BT$ which are cyclic groups of prime order $p$.
If $\KK$ is a normal subloop, then $\KK$ is permutable, and so $\HH\KK$ is a subloop. But
$p < | \HH\KK | < \infty$, a contradiction. Therefore the only normal subloops of $\BT$ are
trivial, and so $\BT$ is a simple loop.
\end{proof}

\begin{corollary}
\thmlabel{trivs}
If $\BT$ is a Bruck-Tarski monster, then $N_r(\BT) = N_m(\BT) = C(\BT) = 1$.
\end{corollary}

\begin{proof}
The right and middle nuclei coincide in any Bol loop and are normal.
The triviality of the commutant follows from Corollary \corref{bruck-comm}.
\end{proof}

Now we will show that Bruck-Tarski monsters exist. Firstly, we need the following.

\begin{lemma}
\lemlabel{invo}
If $\TT$ is a Tarski monster group, then $\Aut(\TT)$ contains no involutions.
\end{lemma}

\begin{proof}
Suppose $\iota \in \Aut(\TT)$ is an involution. Let $G=\TT\langle i\rangle$
be the semidirect product of $\TT$ by $\langle i\rangle$.
Then by (\cite{Ro1}, 14.3.8), either $Z(G)$ contains an involution or
$G$ has a proper infinite subgroup with nontrivial center. In either case,
we have a contradiction.
\end{proof}

\begin{theorem}
\thmlabel{btobt}
If $\TT$ is a Tarski monster, then $\TT(1/2)$ is a Bruck-Tarski monster.
\end{theorem}

\begin{proof}
Suppose that $\KK \neq 1$ is a subloop of $\TT(1/2)$ which is not a cyclic
group of prime order $p$. Then $\KK$ is a twisted subgroup of $\TT$ and
$\TT = \langle\KK\rangle$.
Since $\Aut(\TT)$ contains no involutions (Lemma \lemref{invo}), it follows from
Proposition \prpref{rad} that $\KK' > 1$. But $\TT$ is a simple group, so
$\TT = \KK' \subseteq \KK$. Thus as a loop, $\KK=\TT(1/2)$.
\end{proof}

In 1902, Burnside introduced what is now known as the Burnside Problem \cite{Bu}:
Is a finitely generated periodic group of bounded exponent necessarily finite?
Of course, one can formulate the problem more generally for power-associative
loops as follows:

\smallskip

\emph{Given a class $\PP$ of power-associative loops, is a loop from $\PP$ which
is finitely generated, periodic and of bounded exponent necessarily finite?}

\smallskip

For the class $\PP$ of all groups, the question was answered negatively in 1968
by Adian and Novikov \cite{AN}, and is now also seen to be negative by the existence of Tarski
monsters. On the other hand, the answer is affirmative when $\PP$ is the class of solvable
groups (\cite{Ro1}, 5.4.11) or when $\PP$ is the class of commutative Moufang loops (\cite{Br}, Thm. 11.3). We now have an answer for the class of Bruck loops.

\begin{corollary}
\corlabel{burnside}
The Burnside problem for Bruck loops has a negative answer.
\end{corollary}

\begin{proof}
A Bruck-Tarski monster is a finitely generated, but infinite periodic Bruck loop of bounded exponent.
\end{proof}

We conclude with a few properties of Bruck-Tarski monsters constructed from
Tarski monster groups.
Given a group $G$ and $x,y\in G$, we use the usual notation $x^y = y^{-1}xy$ and
$[x,y] = x^{-1}y^{-1}xy=x^{-1}x^y$.

\begin{theorem}
\thmlabel{bt-comm1}
Let $\TT$ is a Tarski monster. If $1\neq x\in\TT$, then $C_{\TT(1/2)}(x) = \langle x\rangle$.
\end{theorem}

\begin{proof}
Clearly $\langle x\rangle \subseteq C_{\TT(1/2)}(x)$. For the converse, suppose
that $x,y\in \TT(1/2)$ lie in distinct cyclic subloops, and assume that
$y\in C_{\TT(1/2)}(x)$. Then
$x\odot y = y\odot x$ iff $(yx^2y)^{1/2} = (xy^2x)^{1/2}$ iff $yx^2 y = xy^2 x$
iff $[x,y] = [x\inv,y\inv]$. Thus $[x,y]xy = [x\inv,y\inv]xy = xy[x,y]$, that is,
$xy \in C_{\TT}([x,y]) = \langle [x,y]\rangle$. Thus $xy$ generates $\langle [x,y]\rangle$.
Now $(xy)^x = (xy)^{y\inv} = yx = xy [y, x] = xy [x,y]\inv \in \langle [x,y]\rangle$.
Since $x$ and $y$ generate $\TT$, this shows $\langle [x,y]\rangle$ is normal in $\TT$.
Since $\TT$ is simple, either $\langle [x,y]\rangle = 1$ or
$\langle [x,y]\rangle = \TT$, which in either case is a contradiction.
\end{proof}

\begin{theorem}
\thmlabel{bt-lnuc}
If $\TT$ is a Tarski monster, then $N_l(\TT(1/2)) = 1$.
\end{theorem}

\begin{proof}
Assume that $1\neq a\in N_l(\TT(1/2))$. Then for all $x,y\in\TT(1/2)$, we have
$a\odot(x\odot y) = (a\odot x)\odot y$, that is,
$((yx^2y)^{1/2}a^2(yx^2y)^{1/2})^{1/2} = (yxa^2xy)^{1/2}$. Squaring both sides,
we get $(yx^2y)^{1/2}a^2(yx^2y)^{1/2} = yxa^2xy$, or equivalently,
\[
x\inv y\inv (yx^2y)^{1/2} a^2 (yx^2y)^{1/2} y\inv x\inv = a^2\,.
\]
But $x\inv y\inv (yx^2 y)^{1/2} (yx^2y)^{1/2} y\inv x\inv = 1$, and therefore
$x\inv y\inv (yx^2 y)^{1/2} \in C_{\TT}(a^2) = \langle a \rangle$ for all
$x,y\in \TT$.
Suppose that for some $x,y\in \TT$, $x\inv y\inv (yx^2 y)^{1/2}$ generates
$\langle a\rangle$. For every $z\in \TT$,
$(x\inv y\inv (yx^2 y)^{1/2})^z = (x^z)\inv (y^z)\inv (y^z (x^z)^2 y^z)^{1/2}
\in \langle a\rangle$, and so $\langle a\rangle$ is a normal subgroup, a
contradiction. Thus for all $x,y\in \TT$, $x\inv y\inv (yx^2 y)^{1/2} = 1$,
or $yx^2 y = (yx)^2$ or $xy = yx$, another contradiction. Therefore
$N_l(\TT(1/2)) = 1$.
\end{proof}

\textbf{Acknowledgment.}
The authors would like to thank Enrico Jabara and Victor Mazurov
for suggesting the proof of Lemma \lemref{invo}, and Keith Kearnes for
helping us fix Theorem \thmref{bt-comm1}.

\end{document}